\documentclass[11pt]{amsart}

\usepackage{amscd}
\usepackage{amsmath}
\usepackage{graphicx}
\usepackage{amsfonts}
\usepackage{amssymb}
\begin{document}
\title[Generalized Jordan Higher Derivations]{
Characterization of Generalized Jordan Higher Derivations on
Triangular rings}

\author{Xiaofei Qi}
\address[Xiaofei Qi]{
Department of Mathematics, Shanxi University , Taiyuan 030006, P. R.
 China;} \email{qixf1980@126.com}


\thanks{{\it 2010 Mathematical Subject Classification.} 47L35; 47D45}
\thanks{{\it Key words and phrases.} Triangular rings, generalized Jordan higher derivations,
Jordan higher derivation.}
\thanks{This work is supported by National Natural Science Foundation of
China (10771157),  Tianyuan Founds of China (11026161) and
Foundation of Shanxi University.}

\begin{abstract}

Let $\mathcal A$ and $\mathcal B$ be unital rings  and $\mathcal
M$ be a $(\mathcal A, \mathcal B)$-bimodule, which is faithful as
a left $\mathcal A$-module and also as a right $\mathcal
B$-module. Let ${\mathcal U}=\mbox{\rm Tri}(\mathcal A, \mathcal
M, \mathcal B)$ be the associated triangular ring. It is shown
that every additive generalized Jordan (triple) higher derivation
on $\mathcal U$  is a generalized higher derivation.
\end{abstract}
\maketitle

\section{Introduction}
Let $\mathcal A$ be a ring (or an algebra over a commutative ring)
and $\mathcal M$ be an $\mathcal A$-bimodule. Recall that an
additive (linear) map $\delta$ from $\mathcal A$ into $\mathcal M$
is called a derivation if $\delta(AB)=\delta(A)B+A\delta(B)$ for
all $A$,$B\in {\mathcal A}$; a Jordan derivation if
$\delta(A^2)=\delta(A)A+A\delta(A)$ for each $A\in {\mathcal A }$;
and a  Jordan triple derivation if
$\delta(ABA)=\delta(A)BA+A\delta(B)A+AB\delta(A)$ for all
$A$,$B\in {\mathcal A}$. More generally, if there is a  derivation
$\tau:{\mathcal A }\rightarrow {\mathcal M}$ such that
$\delta(AB)=\delta(A)B+A\tau(B)$ for all $A$,$B\in {\mathcal A }$,
then $\delta$ is called a  generalized derivation and $\tau$ is
the relating derivation;  if there is a  Jordan derivation
$\tau:{\mathcal A }\rightarrow {\mathcal M}$ such that
$\delta(A^2)=\delta(A)A+A\tau(A)$ for all $A\in {\mathcal A} $,
then $\delta$ is called a  generalized Jordan derivation and
$\tau$ is the relating Jordan derivation; if there is a  Jordan
triple derivation $\tau:{\mathcal A} \rightarrow {\mathcal M}$
such that $\delta(ABA)=\delta(A)BA+A\tau(B)A+AB\tau(A)$ for all
$A,B\in {\mathcal A }$, then $\delta$ is called a  generalized
Jordan triple  derivation and $\tau$ is the relating Jordan triple
derivation.

 The structures of
derivations, Jordan derivations, generalized derivations and
generalized Jordan derivations were  systematically studied. It is
obvious that every generalized derivation is a generalized Jordan
derivation. But the converse is in general not true.  Zhu in
\cite{ZH2} proved that every generalized Jordan derivation from a
2-torsion free semiprime ring with identity into itself is a
generalized derivation. Hou and Qi in \cite{HQ1} proved that every
additive generalized Jordan derivation of nest algebras on a Banach
space  is an additive generalized derivation. For other results, see
\cite{B2,B3,H2,Lu} and the references therein.

On the other hand, higher derivations had been studied. We first
recall the concepts about higher derivations and generalized
higher derivations.

{\bf Definition 1.1.} (\cite{F2}) {\it Let $D=(\tau_{i})_{i\in
\mathbb{N}}$ be a family of additive maps of ring $\mathcal R$ such
that $\tau_{0}={\rm id}_{\mathcal R}.$ $D$ is said to be:

a higher derivation $( HD,$ for short $)$ if for every $n\in
\mathbb{N}$ we have $\tau_{n}(AB)=\sum
_{i+j=n}\tau_{i}(A)\tau_{j}(B)$ for all $A,B\in \mathcal R$;

a Jordan higher derivation $( JHD,$ for short $)$ if for every $n\in
\mathbb{N}$ we have $\tau_{n}(A^2)=\sum
_{i+j=n}\tau_{i}(A)\tau_{j}(A)$ for all $A\in \mathcal R$;

a Jordan triple higher derivation $( JTHD,$ for short $)$ if for
every $n\in \mathbb{N}$ we have $\tau_{n}(ABA)= \sum
_{i+j+k=n}\tau_{i}(A)\tau_{j}(B)\tau_{k}(A)$ for all $A,B\in
\mathcal R$.}

{\bf Definition 1.2.} (\cite{N}) {\it Let $G=(\delta_{i})_{i\in
\mathbb{N}}$ be a family of additive maps of ring $\mathcal R$ such
that $\delta_{0}={\rm id}_{\mathcal R}.$ $G$ is said to be:

a generalized higher derivation $( GHD,$ for short $)$ if there
exists a higher derivation $D=(\tau_{i})_{i\in \mathbb{N}}$ such
that for every $n\in \mathbb{N}$ we have $\delta_{n}(AB)= \sum
_{i+j=n}\delta_{i}(A)\tau_{j}(B)$ for all $A,B\in \mathcal R$;

a generalized Jordan higher derivation $( GJHD,$ for short $)$ if
there exists a Jordan higher derivation $D=(\tau_{i})_{i\in
\mathbb{N}}$ such that for every $n\in \mathbb{N}$ we have
$\delta_{n}(A^2)= \sum _{i+j=n}\delta_{i}(A)\tau_{j}(A)$ for all
$A\in \mathcal R$;

a generalized Jordan triple higher derivation $( GJTHD,$ for short
$)$ if there exists a Jordan triple higher derivation
$D=(\tau_{i})_{i\in \mathbb{N}}$ such that for every $n\in
\mathbb{N}$ we have $\delta_{n}(ABA)= \sum
_{i+j+k=n}\delta_{i}(A)\tau_{j}(B)\tau_{k}(A)$ for all $A,B\in
\mathcal R$.}

M. Ferrero and C. Haetinger in \cite{F1} proved that every Jordan
higher derivation of a 2-torsion-free ring is a Jordan triple
higher derivation and every Jordan triple higher derivation in a
2-torsion-free semiprime ring is a higher derivation.   Y. S. Jung
in \cite{J1} proved that every generalized Jordan triple higher
derivation on a 2-torsion-free prime ring is a generalized higher
derivation. Recently, Hou and Qi \cite{HQ2} proved that every
additive Jordan (triple) higher derivation of nest algebras on a
Banach space is a higher derivation. Xiao and Wei in \cite{XW}
proved that every Jordan higher derivation on triangular algebras
is a higher derivation.

In the present paper, we will consider generalized Jordan
derivations on triangular rings. In fact, we show that every
additive generalized Jordan higher derivation on triangular rings is
a generalized higher derivation (Theorem 3.1). By using the result,
we  prove that every generalized Jordan triple higher derivation on
triangular rings is also a generalized higher derivation (Theorem
3.2).

Let $\mathcal A$ and $\mathcal B$ be unital rings (or algebras over
a commutative ring $\mathcal R$), and $\mathcal M$ be a $(\mathcal
A, \mathcal B)$-bimodule, which is faithful as a left $\mathcal
A$-module and also as a right $\mathcal B$-module, that is, for any
$a\in {\mathcal A}$ and $b\in \mathcal B$, $a{\mathcal M}={\mathcal
M}b=\{0\}$ imply $a=0$ and $b=0$. The $\mathcal R$-ring ($\mathcal
R$-algebra)
$${\mathcal U}=\mbox{\rm Tri}(\mathcal A, \mathcal M, \mathcal B)=\{
\left(\begin{array}{cc}
a & m\\
  0 & b
\end{array}\right): a\in \mathcal A, m\in \mathcal M, b\in \mathcal B \}$$
under the usual matrix operations is called a triangular ring
(algebra), and the
idempotent element $P=\left(\begin{array}{cc} I_{\mathcal A} &0\\
0&0
\end{array}\right)$ is called the standard idempotent of  ${\mathcal U}$. Clearly, $I-P=\left(\begin{array}{cc} 0 &0\\
0&I_{\mathcal B}
\end{array}\right)$. Here $I$, $I_\mathcal A$ and $I_\mathcal B$ are units of
${\mathcal U}$, ${\mathcal A}$ and $\mathcal B$, respectively. For
more details for triangular rings (algebras) and its relating
questions, the reader see \cite{Ch,QH} and the references therein.

Throughout this paper,  $\mathbb{N}$ denotes the set of natural
numbers including 0.

\section{Preliminaries}

In this section, we give some preliminaries which are needed in
Section 3.

{\bf Lemma 2.1.} {\it Let $\mathcal A$ and $\mathcal B$ be unital
rings, and $\mathcal M$ be a $(\mathcal A, \mathcal B)$-bimodule,
which is faithful as a left $\mathcal A$-module and also as a right
$\mathcal B$-module. Let ${\mathcal U}=\mbox{\rm Tri}(\mathcal A,
\mathcal M, \mathcal B)$ be the triangular ring. Assume that
$G=(\delta_{i})_{i\in \mathbb{N}}$ is an additive generalized Jordan
higher derivation of ${\mathcal U}$ and $D=(\tau_{i})_{i\in
\mathbb{N}}$ the relating additive Jordan higher derivation.  Then
for all $X,Y\in {\mathcal U}$, the following statements hold:}

(1) {\it
$\delta_{n}(XY+YX)=\sum_{i+j=n}(\delta_{i}(X)\tau_{j}(Y)+\delta_{i}(Y)\tau_{j}(X)),$}

(2) {\it
$\delta_{n}(XYX)=\sum_{i+j+k=n}\delta_{i}(X)\tau_{j}(Y)\tau_{k}(X).$}

{\bf Proof.} (1) On the one hand, we have
$$\begin{array}{rl}\delta_{n}((X+Y)^2)
=&\sum_{i+j=n}\delta_{i}(X+Y)\tau_{j}(X+Y)\\
=&\sum_{i+j=n}(\delta_{i}(X)\tau_{j}(X)+\delta_{i}(X)\tau_{j}(Y)\\
 &+\delta_{i}(Y)\tau_{j}(X)+\delta_{i}(Y)\tau_{j}(Y)),
\end{array}$$ and on
the other hand,
$$\begin{array}{rl}\delta_{n}((X+Y)^2)=&\delta_{n}(X^2+XY+YX+Y^2) \\
=&\sum _{i+j=n}\delta_{i}(X)\tau_{j}(X)+\delta_{n}(XY+YX)+\sum
_{i+j=n}\delta_{i}(Y)\tau_{j}(Y).\end{array}$$ Comparing the above
two equations, we obtain that
$$\delta_{n}(XY+YX)=\sum_{i+j=n}(\delta_{i}(X)\tau_{j}(Y)+\delta_{i}(Y)\tau_{j}(X)). $$

(2) Let $S=\delta_{n}(X(XY+YX)+(XY+YX)X).$ By \cite{XW},
$D=(\tau_{i})_{i\in \mathbb{N}}$ is a higher derivation. Then,
using (1) and the fact, on the one hand, we have
$$\begin{array}{rl}
S=&\sum_{i+j=n}(\delta_{i}(X)\tau_{j}(XY+YX)+\delta_{i}(XY+YX)\tau_{j}(X))\\
=&\sum_{i+j=n}\sum_{r+s=j}\delta_{i}(X)(\tau_{r}(X)\tau_{s}(Y)+\tau_{r}(Y)\tau_{s}(X))\\
&+\sum_{i+j=n}\sum_{k+l=i}(\delta_{k}(X)\tau_{l}(Y)+\delta_{k}(Y)\tau_{l}(X))\tau_{j}(X)\\
=&\sum_{i+j=n}\sum_{r+s=j}\delta_{i}(X)\tau_{r}(X)\tau_{s}(Y)
+2\sum_{i+j+k=n}\delta_{i}(X)\tau_{j}(Y)\tau_{k}(X)\\
&+\sum_{i+j=n}\sum_{k+l=i}\delta_{k}(Y)\tau_{l}(X)\tau_{j}(X);
\end{array}$$ on the other hand,
$$\begin{array}{rl}S=&\delta_{n}(X^2Y+2XYX+YX^2)\\
=&\sum_{i+j=n}(\delta_{i}(X^2)\tau_{j}(Y)+\delta_{i}(Y)\tau_{j}(X^2))+2\delta_{n}(XYX)\\
=&\sum_{i+j=n}\sum_{r+s=i}(\delta_{r}(X)\tau_{s}(X)\tau_{j}(Y)\\
&+\sum_{i+j=n}\sum_{k+l=j}\delta_{i}(Y)\tau_{k}(X)\tau_{l}(X)+2\delta_{n}(XYX).
\end{array}$$  These two equations imply that
(2) is true, completing the proof of the lemma. \hfill$\Box$

Now, let $P$ be the standard idempotent of $\mathcal U$.  For the
convenience, in the sequel, let $Q=I-P$. Then ${\mathcal
U}=P{\mathcal U}P+P{\mathcal U}Q+Q{\mathcal U}Q$.

By \cite{XW}, every Jordan higher derivation $D=(\tau_{i})_{i\in
\mathbb{N}}$ on the triangular ring $\mathcal U$ is in fact a higher
derivation and satisfies that
$$\tau_{n}(I)=0\ \ {\rm and} \ \ \tau_{n}(P),\tau_{n}(Q)\in P{\mathcal U}Q\eqno(2.1)$$ for all $n\in
\mathbb{N}$. By the definition of higher derivations, we have
$$\tau_n(XY)=\sum_{i+j=n}\tau_i(X)\tau_j(Y)\quad{\rm for\ \ all }\ \ X,Y\in{\mathcal U}.\eqno(2.2)$$
Thus, for any $X\in{\mathcal U}$, by Eq.(2.1) and noting that
$Q{\mathcal U}P=\{0\}$, we get
$$\begin{array}{rl}\tau_n(PXQ)=&\sum_{i+j=n}\tau_i(P)\tau_j(PXQ)\\
=&P\tau_n(PXQ)+\tau_1(P)\tau_{n-1}(PXQ)+\cdots+\tau_n(P)PXQ\in
P{\mathcal U}Q;\end{array}\eqno(2.3)$$
$$\begin{array}{rl}\tau_n(PXP)=&\sum_{i+j=n}\tau_i(P)\tau_j(PXP)\\
=&P\tau_n(PXP)+\tau_1(P)\tau_{n-1}(PXP)+\cdots+\tau_n(P)PXP\in
P{\mathcal U}P+P{\mathcal U}Q\end{array}\eqno(2.4)$$and
$$\begin{array}{rl}\tau_n(QXQ)=&\sum_{i+j=n}\tau_i(Q)\tau_j(QXQ)\\
=&Q\tau_n(QXQ)+\tau_1(Q)\tau_{n-1}(QXQ)+\cdots+\tau_n(Q)QXQ\in
P{\mathcal U}P+Q{\mathcal U}Q.\end{array}\eqno(2.5)$$

{\bf Remark 2.2.}  By the above analysis, for any Jordan higher
derivation $D=(\tau_{i})_{i\in \mathbb{N}}$ on $\mathcal U$, we
have the following  properties:

${\bf P_D:}$ ({i})$\tau_{n}(I)=0;$ ({ii})$\tau_{n}(P)\in
P{\mathcal U}Q$; ({iii})$\tau_{n}(Q)\in P{\mathcal U}Q$; ({iv})
$\tau_{n}(P{\mathcal U}Q)\subseteq P{\mathcal U}Q$;
  ({v})$\tau_{n}(P{\mathcal U}P)\subseteq P{\mathcal U}P+P{\mathcal
U}Q$ and $\tau_{n}(Q{\mathcal U}Q)\subseteq Q{\mathcal
U}Q+P{\mathcal U}Q$ for each $n\in\mathbb N$.

For any generalized higher derivation $D=(\delta_{i})_{i\in
\mathbb{N}}$ on a triangular ring $\mathcal U$, by the definition,
it is clear that $\delta_{1}$ is a generalized Jordan derivation and
$\tau_1$ the relating Jordan derivation. Hence $\delta_{1}$ is a
generalized derivation by \cite{YX}, that is,
$\delta_1(XY)=\delta_1(X)Y+X\tau_1(Y)$ for $\forall X,Y$, and
satisfies
$$\delta_1(P)\in P{\mathcal
U}P+P{\mathcal U}Q\quad{\rm and }\quad \delta_1(Q)\in P{\mathcal
U}Q+Q{\mathcal U}Q.\eqno(2.6)$$ Thus, by Eq.(2.6) and $\bf P_D$
for $n=1$, we have
$$\delta_1(PXQ)=\delta_1(P)PXQ+P\tau_1(PXQ)\in
P{\mathcal U}Q;$$
$$\delta_1(PXP)=\delta_1(P)PXP+P\tau_1(PXP)\in
P{\mathcal U}P+P{\mathcal U}Q$$and
$$\delta_1(QXQ)=\delta_1(Q)QXQ+Q\tau_1(QXQ)\in
P{\mathcal U}P+Q{\mathcal U}Q.$$

{\bf Remark 2.3.}  By the above argument, for any generalized Jordan
higher derivation $D=(\delta_{i})_{i\in \mathbb{N}}$ on a triangular
ring $\mathcal U$, $\delta_{1}$ is in fact a generalized derivation
and also satisfies the following properties:

${\bf P_1:}$ ({i}) $\delta_1(P) \in (P{\mathcal U}P+P{\mathcal
U}Q)$; ({ii}) $\delta_1(P{\mathcal U}Q)\subseteq P{\mathcal U}Q$;
({iii}) $\delta_1(P{\mathcal U}P)\subseteq P{\mathcal U}P+P{\mathcal
U}Q$ and $\delta_1(Q{\mathcal U}Q)\subseteq Q{\mathcal
U}Q+P{\mathcal U}Q$; ({iv}) $\delta_1(XY)=\delta_1(X)Y+X\tau_1(Y)$
for  $\forall X,Y$.

\section{Characterizations of generalized Jordan Higher Derivations}

In this section, we discuss the generalized Jordan higher
derivations on triangular rings. The following is our main result.

{\bf Theorem 3.1.} {\it  Let $\mathcal A$ and $\mathcal B$ be unital
rings, and $\mathcal M$ be a $(\mathcal A, \mathcal B)$-bimodule,
which is faithful as a left $\mathcal A$-module and also as a right
$\mathcal B$-module. Let ${\mathcal U}=\mbox{\rm Tri}(\mathcal A,
\mathcal M, \mathcal B)$ be the triangular ring and $P$ be the
standard idempotent of $\mathcal U$. Assume that
$G=(\delta_{i})_{i\in \mathbb{N}}$ an additive generalized Jordan
higher derivation of ${\mathcal U}$ and $D=(\tau_{i})_{i\in
\mathbb{N}}$ the relating Jordan higher derivation. Then for any
$X,Y\in {\mathcal U}$ and any $n\in\mathbb N$,  we have
$\delta_n(XY)=\Sigma_{i+j=n}\delta_{i}(X)\tau_{j}(Y)$, that is,
$G=(\delta_{i})_{i\in \mathbb{N}}$ is a generalized higher
derivation.}

{\bf Proof.} We proceed by induction on $n\in \mathbb{N}.$ Assume
that $G=(\delta_{i})_{i\in \mathbb{N}}$ be a generalized Jordan
higher derivation of ${\mathcal U}$ and $D=(\tau_{i})_{i\in
\mathbb{N}}$ the relating Jordan higher derivation.

If $n=1$, by Remark 2.3, $\delta_{1}$ is a generalized derivation
satisfying $\bf P_1$. So the theorem is true in this case.

Now suppose that for any $X, Y\in {\mathcal U}$ and any $m < n,$
$\delta_{m}$ satisfies the following properties:

${\bf P_{m}:}$ ({i}) $\delta_m(P) \in (P{\mathcal U}P+P{\mathcal
U}Q)$; ({ii}) $\delta_m(P{\mathcal U}Q)\subseteq P{\mathcal U}Q$;
({iii}) $\delta_m(P{\mathcal U}P)\subseteq P{\mathcal U}P+P{\mathcal
U}Q$ and $\delta_m(Q{\mathcal U}Q)\subseteq Q{\mathcal
U}Q+P{\mathcal U}Q$; ({iv})
$\delta_{m}(XY)=\sum_{i+j=m}\delta_{i}(X)\tau_{j}(Y)$ for $\forall
X,Y$.

Our aim is to show that $\delta_{n}$ satisfies the following
properties:

${\bf P_{n}:}$ ({i})$\delta_n(P) \in (P{\mathcal U}P+P{\mathcal
U}Q)$; ({ii}) $\delta_n(P{\mathcal U}Q)\subseteq P{\mathcal U}Q$;
({iii}) $\delta_n(P{\mathcal U}P)\subseteq P{\mathcal U}P+P{\mathcal
U}Q$ and $\delta_n(Q{\mathcal U}Q)\subseteq Q{\mathcal
U}Q+P{\mathcal U}Q$; ({iv})
$\delta_n(XY)=\sum_{i+j=n}\delta_{i}(X)\tau_{j}(Y)$ for $\forall
X,Y$.\\ And therefore, $G=(\delta_{i})_{i\in \mathbb{N}}$ be a
generalized higher derivation of $\mathcal U$.  We will prove it by
several steps.

{\bf Step 1.} $\delta_n(P) \in (P{\mathcal U}P+P{\mathcal U}Q)$.

In fact, since $\delta_{i}(P) \in (P{\mathcal U}P+P{\mathcal U}Q)$
for $i=1, 2,..., n-1$ and $\tau_{i}(P) \in P{\mathcal U}Q$ for $i=1,
2,..., n$, we have $\delta_{i}(P)\tau_{j}(P)\in P{\mathcal U}Q $ for
$i=1, 2,..., n-1$ with $i+j=n$ and $P\tau_{n}(P)\in P{\mathcal U}Q$.
Hence $$\delta_{n}(P)=\sum _{i+j=n}\delta_{i}(P)\tau_{j}(P)\in
\delta_{n}(P)P+ P{\mathcal U}Q,$$which implies that $\delta_n(P) \in
(P{\mathcal U}P+P{\mathcal U}Q)$.

{\bf Step 2.} $\delta_{n}(P{\mathcal U}Q)\subseteq P{\mathcal U}Q$.

Take any  $X\in {\mathcal U}$. By Lemma 2.1(1), we have
$$\begin{array}{rl}
\delta_{n}(PXQ)
=&\delta_{n}(PPXQ+PXQP)\\
=&\sum_{i+j=n}(\delta_{i}(P)\tau_{j}(PXQ)+\delta_{i}(PXQ)\tau_{j}(P))     \\
=&\delta_{n}(P)PXQ+\delta_{n}(PXQ)P+P\tau_{n}(PXQ)+PXQ\tau_{n}(P)\\
 &+ \sum_{i+j=n;i\neq
0,n}(\delta_{i}(P)\tau_{j}(PXQ)+\delta_{i}(PXQ)\tau_{j}(P)) .
\end{array}$$
With Step 1 and the properties $\bf P_{m}$, $\bf P_{D}$, it is
clear that $\sum_{i+j=n;i\neq
0,n}(\delta_{i}(P)\tau_{j}(PXQ)+\delta_{i}(PXQ)\tau_{j}(P))\in
P{\mathcal U}Q$, $\delta_{n}(P)PXQ\in P{\mathcal U}Q$,
$P\tau_{n}(PXQ)\in P{\mathcal U}Q$ and $PXQ\tau_{n}(P)=0.$ So we
get $\delta_{n}(PXQ)-\delta_{n}(PXQ)P\in P{\mathcal U}Q$, which
implies that $Q\delta_{n}(PXQ)Q=0$.

Similarly, using the equation $\delta_{n}(PXQ)
=\delta_{n}(QPXQ+PXQQ),$ one can get $P\delta_{n}(PXQ)P=0$.  So
$\delta_{n}(PXQ)=P\delta_{n}(PXQ)Q\in P{\mathcal U}Q$.

{\bf Step 3.} $\delta_{n}(P{\mathcal U}P)\subseteq P{\mathcal
U}P+P{\mathcal U}Q$.

For any $X\in \mathcal U$, by Lemma 2.1(2), we have
$$\begin{array}{rl}
\delta_{n}(PXP)
=\sum_{i+j+k=n}\delta_{i}(P)\tau_{j}(PXP)\tau_{k}(P).
\end{array}$$
By Step 1 and  $\bf P_{m}$, $\bf P_{D}$, for any
$i,j,k\in\{0,1,2,\cdots,n\}$, we  have $\delta_{i}(P)\in P{\mathcal
U}P+P{\mathcal U}Q$ and $\tau_{j}(PXP)\in P{\mathcal U}P+P{\mathcal
U}Q$. It follows that
$$\delta_{i}(P)\tau_{j}(PXP)\tau_{k}(P)\in  (P{\mathcal U}P+P{\mathcal
U}Q)(P{\mathcal U}P+P{\mathcal U}Q)= P{\mathcal U}P+P{\mathcal
U}Q,$$ and so $\delta_{n}(PXP)\in P{\mathcal U}P+P{\mathcal U}Q$.

By a similar argument to that of Step 3, one can check that

{\bf Step 4.} $\delta_{n}(Q{\mathcal U}Q)\subseteq P{\mathcal
U}Q+Q{\mathcal U}Q$.

{\bf Step 5.} For any $X$, $Y\in {\mathcal U}$, the following five
equations hold:

(1) $\delta_{n}(PXPYP)Q=\sum_{i+j=n;i\neq
n}(\delta_{i}(PXP)\tau_{j}(PYP))Q$;

(2) $\delta_{n}(PXPYQ)=\sum_{i+j=n}\delta_{i}(PXP)\tau_{j}(PYQ);$

(3) $\delta_{n}(PXQYQ)=\sum_{i+j=n}\delta_{i}(PXQ)\tau_{j}(QYQ);$

(4) $\delta_{n}(QXQYQ)=\sum_{i+j=n}\delta_{i}(QXQ)\tau_{j}(QYQ)$;

(5) $\sum_{i+j=n}\delta_{i}(PXP)\tau_{j}(QYQ)=0$.

In fact, for any $X,Y\in {\mathcal U}$, by Step 2 and $\bf P_{D},$
we have
$$\begin{array}{rl}
\delta_n(PXPYQ)
=&\delta_n(PXPYQ+PYQPXP)\\
=&\sum_{i+j=n}(\delta_{i}(PXP)\tau_{j}(PYQ)+\delta_{i}(PYQ)\tau_{j}(PXP))\\
=&\sum_{i+j=n}\delta_{i}(PXP)\tau_{j}(PYQ).\end{array}$$ That is,
(2) holds.

Similarly, one can check that (3) and (5) is true.

For (1), we have
$$\begin{array}{rl}
\delta_n(PXP)=&\sum_{i+j+k=n}\delta_{i}(P)\tau_{j}(PXP)\tau_{k}(P)\\
=&\sum_{i+j=n}\delta_{i}(P)\tau_{j}(PXP)P+\sum_{i+j=n-1}\delta_{i}(P)\tau_{j}(PXP)\tau_1(P)\\
&+\sum_{i+j=n-2}\delta_{i}(P)\tau_{j}(PXP)\tau_2(P)+\ldots\\&+
\sum_{i+j=1}\delta_{i}(P)\tau_{j}(PXP)\tau_{n-1}(P)+PXP\tau_{n}(P).\end{array}$$
By induction, the above equation becomes
$$\begin{array}{rl}
\delta_n(PXP)=&\sum_{i+j=n}\delta_{i}(P)\tau_{j}(PXP)P+\delta_{n-1}(PXP)\tau_1(P)\\
&+\delta_{n-2}(PXP)\tau_2(P)+\ldots+\delta_{1}(PXP)\tau_{n-1}(P)+PXP\tau_{n}(P),\end{array}$$
and so
$$\begin{array}{rl}
\delta_n(PXP)Q=&\delta_{n-1}(PXP)\tau_1(P)Q+\delta_{n-2}(PXP)\tau_2(P)Q\\
&+\ldots+\delta_{1}(PXP)\tau_{n-1}(P)Q+PXP\tau_{n}(P)Q.\end{array}$$
Thus for any $X,Y\in {\mathcal U}$, we get
$$\begin{array}{rl}
&\delta_n(PXPYP)Q\\=&\delta_{n-1}(PXPYP)\tau_1(P)Q+\delta_{n-2}(PXPYP)\tau_2(P)Q\\
&+\ldots+\delta_{1}(PXPYP)\tau_{n-1}(P)Q+PXPYP\tau_{n}(P)Q\\
=&\sum_{i+j=n-1}\delta_{i}(PXP)\tau_{j}(PYP)\tau_1(P)Q\\&+\sum_{i+j=n-2}\delta_{i}(PXP)\tau_{j}(PYP)\tau_2(P)Q\\
&+\ldots+\sum_{i+j=1}\delta_{i}(PXP)\tau_{j}(PYP)\tau_{n-1}(P)Q+PXPYP\tau_{n}(P)Q\\
=&PXP(\tau_{n-1}(PYP)\tau_1(P)+\tau_{n-2}(PYP)\tau_2(P)+\ldots+\tau_1(PYP)\tau_{n-1}(P)+PYP\tau_{n}(P))Q\\
&+\delta_1(PXP)(\tau_{n-2}(PYP)\tau_1(P)+\tau_{n-3}(PYP)\tau_2(P)+\ldots+\tau_1(PYP)\tau_{n-2}(P))Q\\
&+\ldots+\delta_{n-2}(PXP)(\tau_1(PYP)\tau_1(P)+PYP\tau_2(P))Q+\delta_{n-1}(PXP)(PYP\tau_1(P))Q\\
=&PXP\tau_{n}(PYP)Q+\delta_1(PXP)\tau_{n-1}(PYP)Q\\
&+\ldots+\delta_{n-2}(PXP)\tau_2(PYP)Q+\delta_{n-1}(PXP)\tau_1(PYP)Q\\
=&\sum_{i+j=n;i\not=n}(\delta_{i}(PXP)\tau_{j}(PYP))Q.\end{array}$$
Hence (1) holds.

Finally, we prove (4). For any $X, Y\in{\mathcal U}$, by Lemma
2.1(2), we have
$$\begin{array}{rl}\delta_n(QXQ)=&\sum_{i+j+k=n}\delta_{i}(Q)\tau_{j}(QXQ)\tau_{k}(Q)\\
=&\sum_{i+j=n}\delta_{i}(Q)\tau_{j}(QXQ)Q+\sum_{i+j+k=n;k\not=0}\delta_{i}(Q)\tau_{j}(QXQ)\tau_{k}(Q).\end{array}$$
Note that, by Steps 1, 4  and the properties of $\tau_i$ (Remark
2.2),
$$\sum_{i+j+k=n;k\not=0}\delta_{i}(Q)\tau_{j}(QXQ)\tau_{k}(Q)\in
(P{\mathcal U}Q+Q{\mathcal U}Q)(P{\mathcal U}Q+Q{\mathcal
U}Q)(P{\mathcal U}Q)=\{0\}.$$ Thus we get
$$\delta_n(QXQ)=\sum_{i+j=n}\delta_{i}(Q)\tau_{j}(QXQ)Q\quad{\rm for \ \ all} \ \ X\in{\mathcal U}.$$
Since $\tau_i(QYQ)=\tau_i(QYQ)Q$, the above equation yields
$$\begin{array}{rl}
\delta_n(QXQYQ)=&\sum_{i+j=n}\delta_{i}(Q)\tau_{j}(QXQYQ)Q\\
=&\sum_{i+j=n}\delta_{i}(Q)(\sum_{p+q=j}\tau_{p}(QXQ)\tau_q(QYQ))Q\\
=&\sum_{i+p+q=n}\delta_{i}(Q)\tau_{p}(QXQ)\tau_q(QYQ)Q\\
=&\sum_{s=0}^{n}[\sum_{i+p=s}(\delta_{i}(Q)\tau_{p}(QXQ)Q)\tau_{n-s}(QYQ)]Q\\
=&\sum_{s=0}^{n}[\delta_s(QXQ)\tau_{n-s}(QYQ)]Q\\
=&\sum_{s=0}^{n}\delta_s(QXQ)\tau_{n-s}(QYQ)=\sum_{i+j=n}\delta_i(QXQ)\tau_j(QYQ).\end{array}$$
It follows that (4) holds.

{\bf Step 6.} $\delta_{n}(XY)=\sum_{i+j=n}\delta_{i}(X)\tau_{j}(Y)$
for all $X,Y\in \mathcal U$, that is, the theorem is true.

We first prove that
$[\delta_{n}(XY)-\sum_{i+j=n}\delta_{i}(X)\tau_{j}(Y)]P=0$. In fact,
for any $X,Y,S\in {\mathcal U}$, by Steps 2-5, on the one hand, we
have
$$\begin{array}{rl}\delta_n(XYPSQ)
=&\delta_n(PXPYPSQ)=\sum_{i+j=n}\delta_i(PXPYP)\tau_j(PSQ)\\
=&\delta_n(PXPYP+PXPYQ+PXQYQ+QXQYQ)PSQ\\
&+\sum_{i+j=n;i\not=n}\delta_i(PXPYP)\tau_j(PSQ)\\
=&\delta_n(XY)PSQ+\sum_{i+j=n;i\not=n}\delta_i(PXPYP)\tau_j(PSQ).
\end{array}$$
On the other hand,
$$\begin{array}{rl}\delta_n(XYPSQ)
=&\delta_n(PXPYPSQ)=\sum_{i+j=n}\delta_i(PXP)\tau_j(PYPSQ)\\
=&\sum_{i+j=n}\delta_i(PXP)\sum_{p+q=j}\tau_p(PYP)\tau_q(PSQ)\\
=&\sum_{i+p+q=n}\delta_i(PXP)\tau_p(PYP)\tau_q(PSQ)\\
=&\sum_{i+p=n}\delta_i(PXP)\tau_p(PYP)PSQ+\sum_{i+p+q=n;q\not=0}\delta_i(PXP)\tau_p(PYP)\tau_q(PSQ)\\
=&\sum_{i+p=n}\delta_i(X)\tau_p(PYP)PSQ+\sum_{s+q=n;q\not=0}\delta_s(PXPYP)\tau_q(PSQ)\\
=&\sum_{i+p=n}\delta_i(X)\tau_p(Y)PSQ+\sum_{s+q=n;q\not=0}\delta_s(PXPYP)\tau_q(PSQ).
 \end{array}$$
The last two equations hold since $\delta_i(PXQ+QXQ)\tau_p(PYP)=0$
 and $\tau_p(PYQ+QYQ)PSQ=0$ for all $i,p$ (by induction on $n$, Steps
2-3 and the property $\bf P_D$). Comparing the above two equations,
we obtain
$[\delta_{n}(XY)-\sum_{i+j=n}\delta_{i}(X)\tau_{j}(Y)]PSQ=0$ for all
$PSQ \in P{\mathcal U}Q$. Since ${\mathcal M}$ is  faithful as a
left $\mathcal A$-module and $Q{\mathcal U}P=\{0\}$, it follows that
$$[\delta_{n}(XY)-\sum_{i+j=n}\delta_{i}(X)\tau_{j}(Y)]P=0. \eqno (3.1)$$

We still need to prove that
$[\delta_{n}(XY)-\sum_{i+j=n}\delta_{i}(X)\tau_{j}(Y)]Q=0.$ For any
$S\in {\mathcal U}$, by Steps 1-5, we have
$$\begin{array}{rl}\delta_n(XYQSQ)
=&\delta_n(PXPYQSQ)+\delta_n(PXQYQSQ)+\delta_n(QXQYQSQ)\\
=&\sum_{i+j=n}\delta_i(PXPYQ)\tau_j(QSQ)+\sum_{i+j=n}\delta_i(PXQYQ)\tau_j(QSQ)\\&
+\sum_{i+j=n}\delta_i(QXQYQ)\tau_j(QSQ)\\
=&\delta_n(PXPYP+PXPYQ+PXQYQ+QXQYQ)QSQ\\&-\delta_n(PXPYP)QSQ+\Delta\\
=&\delta_n(XY)QSQ-\delta_n(PXPYP)QSQ+\Delta,
\end{array}$$
where
$$\begin{array}{rl}
\Delta \doteq
&\sum_{i+j=n;i\not=n}\delta_i(PXPYQ)\tau_j(QSQ)+\sum_{i+j=n;i\not=n}\delta_i(PXQYQ)\tau_j(QSQ)\\&
+\sum_{i+j=n;i\not=n}\delta_i(QXQYQ)\tau_j(QSQ).
\end{array}$$
On the other hand,
$$\begin{array}{rl}&\delta_n(XYQSQ)\\
=&\delta_n(PXPYQSQ)+\delta_n(PXQYQSQ)+\delta_n(QXQYQSQ)\\
=&\sum_{i+j=n}\delta_i(PXP)\tau_j(PYQSQ)+\sum_{i+j=n}\delta_i(PXQ)\tau_j(QYQSQ)\\&
+\sum_{i+j=n}\delta_i(QXQ)\tau_j(QYQSQ)\\
=&\sum_{i+p+q=n}\delta_i(PXP)\tau_p(PYQ)\tau_q(QSQ)
+\sum_{i+p+q=n}\delta_i(PXQ)\tau_p(QYQ)\tau_q(QSQ)\\&
+\sum_{i+p+q=n}\delta_i(QXQ)\tau_p(QYQ)\tau_q(QSQ)\\
=&\sum_{i+p=n}\delta_i(PXP)\tau_p(PYQ)QSQ
+\sum_{i+p=n}\delta_i(PXQ)\tau_p(QYQ)QSQ\\&
+\sum_{i+p=n}\delta_i(QXQ)\tau_p(QYQ)QSQ+\sum_{i+p+q=n;q\not=0}\delta_i(PXP)\tau_p(PYQ)\tau_q(QSQ)\\&
+\sum_{i+p+q=n;q\not=0}\delta_i(PXQ)\tau_p(QYQ)\tau_q(QSQ)
+\sum_{i+p+q=n;q\not=0}\delta_i(QXQ)\tau_p(QYQ)\tau_q(QSQ)\\
=&\sum_{i+p=n}\delta_i(PXP+PXQ+QXQ)\tau_p(PYP+QYQ+PYQ)QSQ\\&
-\sum_{i+p=n}\delta_i(PXP)\tau_p(PYP)QSQ
-\sum_{i+p=n}\delta_i(PXP)\tau_p(QYQ)QSQ\\&
-\sum_{i+p=n}\delta_i(PXQ)\tau_p(PYP)QSQ
-\sum_{i+p=n}\delta_i(PXQ)\tau_p(PYQ)QSQ\\&
-\sum_{i+p=n}\delta_i(QXQ)\tau_p(PYP)QSQ
-\sum_{i+p=n}\delta_i(QXQ)\tau_p(PYQ)QSQ+\Delta\\
=&\sum_{i+p=n}\delta_i(X)\tau_p(Y)QSQ-\sum_{i+p=n}\delta_i(PXP)\tau_p(PYP)QSQ+\Delta\\
=&\sum_{i+p=n}\delta_i(X)\tau_p(Y)QSQ-\delta_n(PXPYP)QSQ+\Delta.
\end{array}$$
Comparing the above two equations, we get
$[\delta_{n}(XY)-\sum_{i+j=n}\delta_{i}(X)\tau_{j}(Y)]QSQ=0$ for all
$S\in {\mathcal U}$. It follows  that
$$[\delta_{n}(XY)-\sum_{i+j=n}\delta_{i}(X)\tau_{j}(Y)]Q=0.\eqno(3.2)$$
Combining Eq.(3.1) and (3.2), we get
$\delta_{n}(XY)=\sum_{i+j=n}\delta_{i}(X)\tau_{j}(Y)$ for all
$X,Y\in{\mathcal U}$. The proof is complete.\hfill$\Box$

Let $F=(\delta_{i})_{i\in \mathbb{N}}$ be any generalized Jordan
triple higher derivation of $\mathcal U$ and $D=(\tau_{i})_{i\in
\mathbb{N}}$ the relating  Jordan triple higher derivation. It is
easy to check  that $\tau_{i}(I)=0 (i=1, 2,...,n)$. So it is obvious
from Lemma 2.1(2) that $F$ is also a generalized Jordan higher
derivation of $\mathcal U$. Hence the following theorem is
immediate.

{\bf Theorem 3.2.} {\it Let $\mathcal A$ and $\mathcal B$ be unital
rings, and $\mathcal M$ be a $(\mathcal A, \mathcal B)$-bimodule,
which is faithful as a left $\mathcal A$-module and also as a right
$\mathcal B$-module. Let ${\mathcal U}=\mbox{\rm Tri}(\mathcal A,
\mathcal M, \mathcal B)$ be the triangular ring. Then every
generalized Jordan triple higher derivation of $\mathcal U$ is a
generalized higher derivation.}

Recall that a nest $\mathcal N$ on a Banach space $X$ is a chain
of closed  subspaces of $X$ which is closed under the formation of
arbitrary closed linear span  and intersection, and which includes
$\{0 \}$ and $X$. The nest algebra associated to the nest
$\mathcal N$, denoted by Alg$\mathcal N$, is the weakly closed
operator algebra consisting of all operators that leave $\mathcal
N$ invariant, i.e.,
$$\mbox{\rm Alg}\mathcal N=\{T\in \mathcal B(X): TN\subseteq N  \mbox{ for all }  N\in
\mathcal N \}.$$  If $X$ is a Hilbert space, then every
$N\in{\mathcal N}$ corresponds to a projection $P_N$ satisfying
$P_N=P_N^*=P_N^2$ and $N=P_N(X)$. However, it is not always the
case for general  nests on Banach spaces as $N\in{\mathcal N}$ may
be not complemented. We refer the reader to \cite{D} for the
theory of nest algebras.

As an application of Theorem 3.1 and 3.2 to the nest algebra case,
we have

{\bf Theorem 3.3.} {\it Let $\mathcal N$ be a nest on a Banach
space $X$ and there exists a non-trivial element in $\mathcal N$
which is complemented in $X$ $($ in particular, $\mathcal N$ be
any nest on a Hilbert space $H$ $)$. Then every generalized Jordan
(triple) higher derivation on $\mbox{\rm Alg}\mathcal N$ is a
generalized higher derivation.}

{\bf Remark 3.4.}  By Theorem 3.1 and  3.2, we show that the
concepts of generalized Jordan triple higher derivation and
generalized Jordan higher derivation on triangular algebras   are
equivalent, and so the concepts of generalized Jordan triple higher
derivation, generalized Jordan higher derivation and generalized
higher derivation on triangular rings are equivalent to each other.

\end{document}